\documentclass[12pt]{article}
\usepackage{amsmath}
\usepackage[cp1251]{inputenc}
\usepackage{amscd}
\usepackage{amssymb}
\begin{document}

\noindent{\large\bf Rate of approximation by logarithmic
derivatives of polynomials whose zeros lie on a
circle}\footnote{This work supported by RFBR project
18-31-00312 mol\underline{\ \ }a.} \\

\noindent{\bf M. A. Komarov}\footnote{Department of Functional
Analysis and Its Applications, Vladimir State University,
Gor$'$kogo street 87,
600000 Vladimir, Russia\\ e-mail: kami9@yandex.ru}\\

\noindent{\bf Abstract} \ We obtain an estimate for uniform
approximation rate of bounded analytic in the unit disk functions
by logarithmic derivatives of $C$-polynomials, i.e., polynomials,
all of whose zeros lie on the unit circle $C:|z|=1$.

\medskip

$C$-polynomials, uniform
approximation

\medskip


{\bf 1.} Let $C$ denote the unit circle $|z|=1$ and $D$ denote the
unit disk $|z|<1$. It's proved in \cite{Thompson} for any bounded
analytic in $D$ function $f$ and in \cite{Rub-Suff} for any
analytic in $D$ function, that there is a sequence of rational
functions $S_n(z)=S_n(f;z)$ of the form
\[S_n(z)=\sum_{k=1}^{m_n}(z-z_{n,k})^{-1}, \qquad |z_{n,k}|=1,\]
which converges to $f(z)$ uniformly on every closed subset of $D$.
Obviously, $S_n$ is a logarithmic derivative of $C$-polynomial
$(z-z_{n,1})\dots(z-z_{n,m_n})$ ($C$-{\it polynomials} defined in
abstract). Analogous problems with more general constraints on
poles (for example, if $z_{n,k}$ belong to rectifiable Jordan
curve) investigated in \cite{Korev,Chui,Chui-Shen}. Approximation
by sums $\sum_{k=1}^n(z-z_k)^{-1}$ with {\it free} poles studied
in \cite{DD,Kos,K-IzvRAN-2017} (see also bibliography in
\cite{K-IzvRAN-2017}).

We study a rate of approximation of bounded analytic in $D$
functions by logarithmic derivatives of $C$-polynomials on closed
subsets of $D$.

\medskip

\noindent{\bf Theorem.} {\it For any bounded analytic in $D$
function $f(z)$ there is a sequence of $C$-po\-ly\-no\-mi\-als
$P_{N}(z)$, $N\ge N_0$, such that ${\rm deg}\,P_N(z)=N$ and
\begin{equation}\label{result}
    \left|\frac{P_N'(z)}{P_N(z)}-f(z)\right|<
    \frac{(a+\varepsilon)^{n+1}}{\varepsilon(1-a-\varepsilon)}
    (1+o(1)), \qquad n=[N/2], \qquad {\rm as} \quad N\to \infty
\end{equation}
in any disk $K_a=\{|z|\le a\}$ $(a<1)$ for every $\varepsilon\in
(0,1-a)$.}

\medskip

Let $d_n(f,K_a)$ denote the error in best approximation to $f$ on
the disk $K_a$ by logarithmic derivatives $Q'/Q$ of polynomials
$Q$ of degree at most $n$ with {\it free} zeros. It's interesting,
that, generally speaking, the convergence $d_n(f,K_a)$ is also
geometric (sf. (\ref{result})):
\[\limsup_{n\to \infty} \sqrt[n]{d_n(f,K_a)}\le a.\]
This estimate follows from \cite{Kos} (see also \cite{DD}), where
the analog of polynomial Walsh's theorem was proved for the
problem of approximation by such fractions $Q'/Q$.

\medskip

{\bf 2.} To construct polynomials $P_N(z)$ we use the approach
\cite{Rub}. We need the next lemma, stated in \cite{Rub} (for the
case $m=0$ see \cite[p.\,108]{Polya-Szego}). Further
$\overline{D}=D+C$.

\medskip

\noindent{\bf Lemma} \cite{Rub}. {\it Let
$Q(z)=a_0(z-z_1)\dots(z-z_q)$, $a_0\ne 0$, be a polynomial of
degree $q\ge 1$ and
$Q^*(z):=z^q\overline{Q}(1/z)=\overline{a}_0(1-\overline{z}_1
z)\dots (1-\overline{z}_q z)$. If $Q(z)$ zero free in
$\overline{D}$, then $Q(z)+z^m Q^*(z)$ is $C$-polynomial for every
$m=0,1,2,\dots$, and $|Q^*(z)|\le |Q(z)|$ in $\overline{D}$.}

\medskip

To prove this lemma it is sufficient to consider the equation
\begin{equation}\label{Phi=...}
    \Phi(z)=-a_0/\overline{a}_0, \qquad
    \Phi(z):=\frac{a_0}{\overline{a}_0}\frac{z^m Q^*(z)}{Q(z)}\equiv
    z^m \prod_{j=1}^q \frac{1-\overline{z}_j z}{z-z_j}, \qquad
    |z_j|>1.
\end{equation}
Absolute values of all factors in last product are less, equal or
more than 1 iff $|z|<1$, $|z|=1$ or $|z|>1$, respectively,
therefore all roots of equation (\ref{Phi=...}) (and so, all zeros
of $Q(z)+z^m Q^*(z)$) lie on $C$, and in the disk $\overline{D}$
we have $|\Phi(z)|\le 1$ and $|Q^*(z)|\le |Q(z)|$.

\medskip

\noindent{\it Remark 1.} If $Q(z)\equiv a_0={\rm const}\ne 0$,
then $Q^*(z)\equiv \overline{a}_0$, consequently, $Q(z)+z^m
Q^*(z)\equiv a_0+z^m \overline{a}_0$ is $C$-polynomial for $m>0$
only.

\smallskip

\noindent{\it Remark 2.} Lemma is also true, if zeros of $Q(z)$
(but not all of them) lie on $C$, because $\overline{z}_j=1/z_j$
and $|1-\overline{z}_j z|/|z-z_j|\equiv 1$ as $|z_j|=1$. But if
$Q(z)$ is $C$-polynomial, then $Q^*(z)\equiv tQ(z)$,
$t=\overline{Q(0)}/a_0$, and we need again to assume $m>0$, if
$1+t=0$. \medskip

{\bf 3.} {\it Proof of the theorem.} Set $g(z)=\exp\left(\int_0^z
f(\zeta)d\zeta\right)$,
\[g(z)=s_n(z)+R_{n}(z), \qquad s_n(z)=1+\sum_{k=1}^n g_k z^k, \qquad
R_n(z)=\sum_{k=n+1}^{\infty} g_k z^k.\] Derivative $g'\equiv gf$
is bounded and analytic in $D$. In particular, $g'$ belongs to the
Hardy space $H^1(D)$, hence the series $\sum |g_k|$ converges
\cite[Theorem\,15]{Hardy}.

Function $\varphi(x,y)={\rm Re} \int_0^z f(\zeta)d\zeta$ is
bounded in $D$, so $\varphi>-\infty$ and $\inf_{D} |g(z)|=\inf_{D}
\exp \varphi(x,y)=M_0>0$. Choose $n_0\in \mathbb{N}$, such that
$\sum_{n+1}^{\infty} |g_k|\le M_0-M_0/2$ as $n\ge n_0$. Hence
polynomials $s_n(z)$ are zero free in $\overline{D}$ as $n\ge
n_0$.

Further $N\ge 2n_0$ and $n:=[N/2]\ge n_0$. Set
$p(z)=s_n^*(z)=z^{q}\overline{s}_n(1/z)$, where $q=\deg s_n(z)$,
$0\le q\le n$. By lemma and remark 1 we have $|p(z)|\le |s_n(z)|$
in $\overline{D}$, and sums
\[P(z)=P_{q+m}(z):=s_n(z)+z^{m}p(z) \qquad {\rm as} \quad m=1,2,\dots\]
are $C$-polynomials. Rewrite $P$ as $P(z)\equiv
g(z)+z^{m}p(z)-R_n(z)$. We now have
\[P'(z)=g(z)f(z)+m z^{m-1}p(z)+z^{m}p'(z)-R_n'(z),\]
\begin{equation}\label{P'/P-f}
    \frac{P'(z)}{P(z)}-f(z)= \frac{z^{m-1}(m
    -zf(z))p(z)+z^{m}p'(z)-R_n'(z)+f(z)R_n(z)}{P(z)}.
\end{equation}
Denote $M_1=\max\{1,|g_1|,\dots,|g_n|\}$. Since $|g_k|<M_0$ (as
$k\ge n_0+1$), we have
\begin{equation}\label{|pn|,|Rn|}
|p(z)|\le |s_n(z)|<M_1/(1-a), \qquad |R_n(z)|<M_0 a^{n+1}/(1-a)
\qquad {\rm as} \quad |z|\le a,
\end{equation}
\begin{equation}\label{|P|>}
    |P(z)|\ge |g(z)|-|z^{m}p(z)-R_n(z)|>
    M_0-(M_1a^m+M_0a^{n+1})/(1-a) \qquad {\rm as} \quad |z|\le a.
\end{equation}
If $|z|<r<1$ and function $F$ is analytic in $D$, then
\[2\pi|F'(z)|=\left|\int_{|\zeta|=r}\frac{F(\zeta)d\zeta}{(\zeta-z)^2}\right|
\le \max_{|\zeta|=r}|F(\zeta)|
\int_{|\zeta|=r}\frac{|d\zeta|}{|\zeta-z|^2}=
\max_{|\zeta|=r}|F(\zeta)|\frac{2\pi r}{r^2-|z|^2}\] (we apply
Poisson's integral), and hence if $r=a+\varepsilon<1$, then
\[|F'(z)|<\varepsilon^{-1}\max\nolimits_{|\zeta|=a+\varepsilon}|F(\zeta)|
\qquad {\rm as} \quad |z|\le a.\] Thus, by this and
(\ref{|pn|,|Rn|}) we have
\begin{equation}\label{|p'|<}
    |p'(z)|<\frac{M_1}{\varepsilon(1-a-\varepsilon)}, \qquad
    |R_n'(z)|<\frac{M_0(a+\varepsilon)^{n+1}}{\varepsilon(1-a-\varepsilon)}
    \qquad {\rm as} \quad |z|\le a.
\end{equation}
We put $m=N-q\ge n$ and obtain (\ref{result}) from
(\ref{P'/P-f})---(\ref{|p'|<}), and theorem follows.

\medskip

\noindent{\it Remark 3.} $P_N'(z)/P_N(z)-f(z)=O(z^l)$ as $l\ge
n-1$ (see (\ref{P'/P-f})).

\smallskip

\noindent{\it Remark 4.} It follows from $g'\in H^1(D)$, that
$g(z)$ continuous in $\overline{D}$ and absolutely continuous on
$C$ \cite[Ch.II, \S5(5.7)]{Priwalow}. In particular,
$g(e^{i\theta})$ has bounded variation and $|g_k|=O(1/k)$, so we
can to write $O(1/n)$ instead of $1+o(1)$ in (1). If $g$ is a zero
free in $\overline{D}$ polynomial of degree $q\ge 0$ and $f=g'/g$,
then $R_n(z)\equiv 0$ as $n\ge q$, and approximation rate is
higher. For example, if $f(z)\equiv 0$, then $P_N(z)=1+z^N$ and
$\sup_{K_a}|{P}_N'/{P}_N-f|=Na^{N-1}/(1-a^N)$.

\end{document}